\documentclass[12pt,a4paper]{article}
\usepackage[english]{babel}
\usepackage{amsmath}
\usepackage{amssymb}
\usepackage{amsfonts}
\usepackage{hyperref}

\begin{document}


\begin{center}
{\huge Modeling viscous compressible barotropic multi-fluid flows}
\end{center}

\medskip

\begin{center}
{\large Alexander Mamontov,\quad Dmitriy Prokudin\footnote{The authors were supported by the Russian Science Foundation (grant 15--11--20019).}}
\end{center}

\medskip

\begin{center}
{\large August 25, 2017}
\end{center}

\medskip

\begin{center}
{
Lavrentyev Institute of Hydrodynamics, \\ Siberian Branch of the Russian Academy of Sciences\\ pr. Lavrent'eva 15, Novosibirsk 630090, Russia}
\end{center}

\medskip

\begin{center}
{\bfseries Abstract}
\end{center}


\begin{center}
\begin{minipage}{110mm}
 We study the system of equations which describes barotropic (isentropic) flows of viscous compressible multi-fluids (mixtures of fluids). We study the relations between pressure, densities, concentrations, viscosities and other parameters of the flow and of the medium.
\end{minipage}
\end{center}

\bigskip


\newpage

\tableofcontents

\bigskip

\section{Introduction}
The flow of a mixture which consists of $N$ viscous compressible fluids (gases) with a common pressure of the constituents (phases) may be described \cite{mamprok.france}, \cite{mamprok.semr17} via the system of equations
\begin{equation}\label{mamprok3d.continit}
\frac{\partial\rho_{i}}{\partial t}+{\rm div\,}(\rho_{i}\boldsymbol{v})=0,\quad i=1, \ldots, N,
\end{equation}
\begin{equation}\label{mamprok3d.mominit}
\frac{\partial(\rho_{i}\boldsymbol{u}_{i})}{\partial t}+{\rm div\,}(\rho_{i}\boldsymbol{v}\otimes\boldsymbol{u}_{i})+\alpha_i\nabla p
={\rm div\,}{\mathbb S}_{i}+\rho_{i}\boldsymbol{f}_{i},\quad i=1, \ldots, N.
\end{equation}
Here $\rho_i$ is the density of the $i$-th component (constituent) of the multi-fluid, $\boldsymbol{u}_{i}$ is the velocity field, and $p$ is the pressure. The concentration
$\xi_i=\rho_i/\rho$ of each component is defined (they compose the vector $\boldsymbol{\xi}=\{\xi_i\}_{i=1}^N$), where $\displaystyle \rho=\sum\limits_{i=1}^N \rho_i$ is the total density of the multi-fluid. Then, $\displaystyle \boldsymbol{v}=\sum\limits_{i=1}^N\alpha_i \boldsymbol{u}_{i}$ is the average velocity of the multi-fluid, and the coefficients $\alpha_i(\boldsymbol{\xi})>0$ are such that $\displaystyle \sum\limits_{i=1}^N \alpha_i=1$ (e.~g., $\alpha_i=\xi_i$). The viscous stress tensors ${\mathbb S}_{i}$ are defined by the equalities
$${\mathbb S}_{i}=\sum\limits_{j=1}^{N}\widehat{{\mathbb S}}_{ij},\quad \widehat{{\mathbb S}}_{ij}=\Big(2\mu_{ij}{\mathbb D}(\boldsymbol{u}_{j})+\lambda_{ij}({\rm div\,}
\boldsymbol{u}_{j}){\mathbb I}\Big),\qquad  i, j=1, \ldots, N,$$
where ${\mathbb D}(\boldsymbol{w})=((\nabla\otimes\boldsymbol{w})+(\nabla\otimes\boldsymbol{w})^{*})/2$ is the rate of deformation tensor of the vector field~$\boldsymbol{w}$,
${\mathbb I}$~is the identity tensor, and the viscosity coefficients compose the matrices
$$\textbf{M}=\{\mu_{ij}\}_{i, j = 1}^{N}>0,\quad \textbf{H}=\boldsymbol{\boldsymbol{\Lambda}}+\frac{2}{3}\textbf{M}\geqslant 0,\quad
\boldsymbol{\Lambda}=\{\lambda_{ij}\}_{i,j= 1}^{N}.$$
Finally, $\boldsymbol{f}_{i}$ are known fields of external body forces.

Generally speaking, the pressure $p$ depends on the densities and on other thermodynamical variables (say, temperature(s)), so that the system \eqref{mamprok3d.continit}, \eqref{mamprok3d.mominit} must be complemented by the equation(s) for these variables (temperatures, energies, entropies etc.). However, the paper is limited by barotropic flows, so we concentrate on the justification of the state (constitutive) equations which characterize the relations between the following values: $\{\rho_i\}_{i=1}^N$, $p$, $\{\alpha_i\}_{i=1}^N$, $\{\lambda_{ij}\}_{i,j=1}^N$ and $\{\mu_{ij}\}_{i,j=1}^N$.

\section{Pressure as a function of the densities}
From (\ref{mamprok3d.continit}) we have
\begin{equation}\label{mamprok3d.contcons}
\frac{\partial\rho}{\partial t}+{\rm div\,}(\rho\boldsymbol{v})=0,\qquad\frac{\partial\xi_{i}}{\partial t}+\boldsymbol{v}\cdot\nabla\xi_{i}=0.
\end{equation}

\subsection{Motion of a material volume}
Let the material volume $V$ be occupied by a mixture of $N$ ideal gases, and this multi-fluid volume is the subject of a quasistatic process, in which the thermodynamical variables are supposed to be constant over the whole volume. Concerning the system \eqref{mamprok3d.continit}, \eqref{mamprok3d.mominit}, this means that the volume is rather small and consists of the constant ensemble of particles, i.~e. it moves along the trajectories of the transport equations \eqref{mamprok3d.continit}.

Each constituent of the multi-fluid is characterized by the following variables (parameters): density $\rho_i$, temperature $\theta$ (which is common for all constituents), entropy $S_i$, mass $m_i=\rho_i V$, molar mass $M_i$, partial pressure $\displaystyle p_i=\frac{\rho_i}{M_i}R\theta$ ($R$ is the gas constant), degree of freedom of molecules $\nu_i$,
adiabatic index (heat capacity ratio, Poisson constant) $\displaystyle \gamma_i=1+\frac{2}{\nu_i}$ (i.~e. $\displaystyle \nu_i=\frac{2}{\gamma_i-1}$), and internal energy $\displaystyle U_i=\frac{\nu_i R}{2M_i}m_i\theta$, or specific internal energy  $\displaystyle e_i=\frac{\nu_i R}{2M_i}\theta$.

The heat transfer into each constituent (from other constituents and from outside of the volume) equals
$$\theta dS_i=dU_i+p_i dV,$$
and the total heat transfer into the volume equals (here $S$ is the total entropy of the volume)
$$\delta Q=\theta dS=\sum\limits_{i=1}^N dU_i+\sum\limits_{i=1}^N p_i dV.$$

\subsection{Adiabatic process}
If the process is adiabatic (isentropic), then $\delta Q=0$, and we obtain
\begin{equation}\label{mamprok3d.first}
\sum\limits_{i=1}^N dU_i+\sum\limits_{i=1}^N p_i dV=0.
\end{equation}
Observe that
\begin{equation}\label{mamprok3d.second}
dU_i+p_i dV=R\theta\left(\frac{m_i \nu_i}{2M_i}d\ln\theta+\frac{m_i}{M_i}d\ln V\right).
\end{equation}
Let us denote
$$\alpha=\sum\limits_{i=1}^N \frac{m_i \nu_i}{2M_i},\qquad \beta=\sum\limits_{i=1}^N \frac{m_i}{M_i}.$$
We find from (\ref{mamprok3d.first}) and (\ref{mamprok3d.second}) that $\theta=C_1 V^{-\beta/\alpha}$. Taking into account the relation $m=\rho V$ ($m$ is the total mass of the multi-fluid), we finally obtain $\theta=C_2 \rho^{\beta/\alpha}$.

\subsection{Adiabatic index}
Let us note that the value
$$\frac{\beta}{\alpha}=\frac{\displaystyle \sum\limits_{i=1}^N \frac{\rho_i}{M_i}}{\displaystyle \sum\limits_{j=1}^N \frac{\rho_j \nu_j}{2M_j}}=
\frac{\displaystyle \sum\limits_{i=1}^N \frac{\xi_i}{M_i}}{\displaystyle \sum\limits_{j=1}^N \frac{\xi_j \nu_j}{2M_j}}$$
is not constant in the flow, however, in the partial case of multi-fluid with equal degree of freedom of molecules $\nu_i=\nu$ in all constituents (so that $\gamma_i=\gamma$) we obtain
$$\frac{\beta}{\alpha}=\frac{2}{\nu}=\gamma-1.$$

In a general case, $\beta/\alpha$ depends on the concentrations in a spacial point and can take any values from the interval $[\gamma_{\rm min}-1,\gamma_{\rm max}-1]$ (where $\gamma_{\rm min}$ and $\gamma_{\rm max}$ are correspondingly the minimal and maximal adiabatic indices between all constituents of the multi-fluid). Let us denote
\begin{equation}\label{mamprok3d.gamma}
\gamma=\gamma(\boldsymbol{\xi})=\frac{\beta}{\alpha}+1=1+\frac{\displaystyle \sum\limits_{i=1}^N \frac{\xi_i}{M_i}}{\displaystyle \sum\limits_{j=1}^N \frac{\xi_j}{M_j(\gamma_j-1)}}.
\end{equation}
This quantity is a function of the concentrations, and hence (see \eqref{mamprok3d.contcons}$_2$) is transported along the trajectories, i.~e. it satisfies the equation
\begin{equation}\label{mamprok3d.transportgamma}
\frac{\partial\gamma}{\partial t}+\boldsymbol{v}\cdot\nabla\gamma=0,
\end{equation}
but $\gamma$ can be supposed to be constant in two situations: either if all $\gamma_i$ coincide (then $\gamma$ coincides with it), or if $\gamma={\rm const}$ in the initial state of the multi-fluid.

\subsection{Completion of the constitutive equation}
Thus, we consider an adiabatic movement of a material volume of a multi-fluid (which consists of the constant ensemble of particles), in which the distribution of the intensive variables is supposed to be uniform. As a result, the partial pressures have taken the form $\displaystyle p_i=C_3 \frac{\rho_i}{M_i}\rho^{\beta/\alpha}$, and hence the total pressure
$$p=\sum\limits_{i=1}^N p_i=C_3 \rho^{\beta/\alpha}\sum\limits_{i=1}^N \frac{\rho_i}{M_i}.$$
Thus,
\begin{equation}\label{mamprok3d.pressure}
p=K_1\rho^{\gamma-1}\widetilde{\rho}=K_1\rho^\gamma\widetilde{\xi},\qquad \widetilde{\rho}=\sum\limits_{i=1}^N \frac{\rho_i}{M_i},
\qquad \widetilde{\xi}=\sum\limits_{i=1}^N \frac{\xi_i}{M_i}.
\end{equation}
Observe that $C_1=V_0^{\beta_0/\alpha_0}\theta_0$, $C_2=\rho_0^{-\beta_0/\alpha_0}\theta_0$, $C_3=K_1=R \rho_0^{-\beta_0/\alpha_0}\theta_0$.

\subsection{Variant of writing the constitutive equation}
Only one of the densities is independent, since they are connected through the relations
$$\rho_i=\frac{m_i}{V}=\frac{m_i}{m}\rho=\frac{\rho_{0i} V_0}{\rho_0 V_0}\rho=\frac{\rho_{0i}}{\rho_0}\rho,$$
and we can choose either one of them as an unknown value, or some their combination, e.~g., the total density $\rho$. As a result, we obtain
$$p_i=\frac{R\theta_0}{M_i}\rho_0^{-\gamma_0}\rho_{0i}\rho^\gamma\quad \Longrightarrow \quad p=K \rho^\gamma, \quad
  K=R\theta_0\rho_0^{-\gamma_0}\sum\limits_{i=1}^N\frac{\rho_{0i}}{M_i}.$$
In this way, we see that, in a real flow, $K$ is a function of the Lagrangian spacial coordinates, i.~e. it is transported along the trajectories and satisfies the equation
\begin{equation}\label{mamprok3d.transportk}
\frac{\partial K}{\partial t}+\boldsymbol{v}\cdot\nabla K=0,
\end{equation}
but, as the simplest version, $K$ can be accepted to be constant (this means a property of the initial distribution of the parameters of the multi-fluid).

\subsection{Choice of the form of the constitutive equation}
The relations (see Section 2.5)
\begin{equation}\label{mamprok3d.pressure2}
p=K \rho^\gamma
\end{equation}
and (see \eqref{mamprok3d.pressure})
\begin{equation}\label{mamprok3d.pressure1}
p=K_1\rho^\gamma\widetilde{\xi},
\end{equation}
are, of course, equivalent (since they are obtained from the same assumptions), as one can easily make sure via the direct verification of the identity $K=K_1\widetilde{\xi}$.
Those relations are just the same, written in two different forms, as long as we discuss an adiabatic movement of a material volume. However, when we use the constitutive equation for the pressure in order to complete the system \eqref{mamprok3d.continit}, \eqref{mamprok3d.mominit}, the choice of the appropriate form may become crucial (see Section 5.4).

\subsection{Coefficients $\alpha_i$}
The quantities $\alpha_i$ which enter as factors at the pressures in \eqref{mamprok3d.mominit}, and also in the relations between $\boldsymbol{v}$ and $\boldsymbol{u}_{i}$, depend on the concentrations, and hence satisfy the equations
\begin{equation}\label{mamprok3d.transportalpha}
\frac{\partial \alpha_i}{\partial t}+\boldsymbol{v}\cdot\nabla\alpha_i=0.
\end{equation}
In the simplest case, $\alpha_i$ can be accepted to be constant (this means a property of the initial state of the multi-fluid).

\section{Viscosities}
The viscosity matrices in real multi-fluids take non-trivial forms. For instance, there is the following information \cite{mamprok.alsharif}, \cite{mamprok.benner}, \cite{mamprok.lee} about the matrix $\textbf{M}$. As~$i\neq j$, the following representation holds
$$\mu_{ij}=\mu_{ij}^0(\widehat{\mu}_i,\widehat{\mu}_j)\xi_i\xi_j\exp\left(\frac{\alpha_{ij}\xi_i+\beta_{ij}\xi_j}{\xi_i+\xi_j}\right),$$
where $\widehat{\mu}_i$ are the ``pure'' viscosities of the constituents. The most frequent variant is
\begin{equation}\label{mamprok3d.coeffsimp}
\mu_{ij}^0(\widehat{\mu}_i,\widehat{\mu}_j)=\sqrt{\widehat{\mu}_i \widehat{\mu}_j},\qquad \alpha_{ij}=\beta_{ij}=0,
\end{equation}
but, generally speaking, $\alpha_{ij}$ and $\beta_{ij}$ are the empiric constants. The diagonal entries of the matrix~$\textbf{M}$ are accepted as
$$\mu_{ii}=\widehat{\mu}_i \xi_i^2+\sum\limits_{j\neq i}\mu_{ij}.$$
Thus, in the simplest case (\ref{mamprok3d.coeffsimp}), we can write
\begin{equation}\label{mamprok3d.viscfin}
\mu_{ij}=\widehat{\nu}_i\widehat{\nu}_j+\widehat{\nu}_i(\widehat{\nu}-\widehat{\nu}_i)\delta_{ij},\qquad \widehat{\nu}_i=\sqrt{\widehat{\mu}_i}\xi_i,\qquad \widehat{\nu}=\sum\limits_{i=1}^N \widehat{\nu}_i.
\end{equation}
The matrix $\textbf{M}$, on the one hand, possesses an important diagonal predominance property, but on the other hand, it essentially depends on the concentrations.

Summarizing, we see that the viscosity matrices, in general, depend on the concentrations, and hence (see \eqref{mamprok3d.contcons}$_2$) are transported along the trajectories and satisfy the equations
\begin{equation}\label{mamprok3d.transportvisc}
\frac{\partial\textbf{M}}{\partial t}+(\boldsymbol{v}\cdot\nabla)\textbf{M}=0,\qquad \frac{\partial\boldsymbol{\Lambda}}{\partial t}+(\boldsymbol{v}\cdot\nabla)\boldsymbol{\Lambda}=0,
\end{equation}
but they can be accepted to be constant as a simplifying assumption which is justified at least in the case of the corresponding initial configuration.

The dependence of the viscosity matrices on the concentrations constitutes serious difficulties for the effective viscous flux (EVF) technique, which underlies the present existence theory for the equations of compressible viscous fluids. That is why, in the present theory, the simplifying assumption of constant viscosities is always accepted, but with the preservation of necessary properties of positiveness.

\section{Energy balance}
Regardless of the variant of the constitutive relation for the pressure (in the framework described in Section 2.6) and of the specification of the values $\{\alpha_i\}_{i=1}^N$, $\{\lambda_{ij}\}_{i,j=1}^N$ and $\{\mu_{ij}\}_{i,j=1}^N$, the total energy balance can be expressed as an integral relation. Indeed, let \eqref{mamprok3d.pressure2} hold, where $K$ and $\gamma$ are arbitrary constants satisfying \eqref{mamprok3d.transportgamma} and \eqref{mamprok3d.transportk}. Observe that $p$ in \eqref{mamprok3d.pressure1} satisfies these assumptions as well, as long as $K_1$ satisfies \eqref{mamprok3d.transportk}, and $\widetilde{\xi}$ satisfies \eqref{mamprok3d.contcons}$_2$. Then, for any smooth function $\varphi$ we have
$$d(\varphi(\gamma)p)=p\varphi'(\gamma)d\gamma+\varphi(\gamma)\rho^\gamma dK+\varphi(\gamma)K\rho^\gamma(\ln\rho)d\gamma+\varphi(\gamma)K\gamma \rho^{\gamma-1}d\rho.$$
Dividing by $dt$, using the notation $\displaystyle \frac{d}{dt}=\frac{\partial}{\partial t}+\boldsymbol{v}\cdot\nabla$ and the equation \eqref{mamprok3d.contcons}$_1$, we obtain
\begin{equation}\label{mamprok3d.transpres}
\frac{\partial}{\partial t}(\varphi(\gamma)p)+{\rm div\,}(\varphi(\gamma)p\boldsymbol{v})+\varphi(\gamma)(\gamma-1)p{\rm div\,}\boldsymbol{v}=0.
\end{equation}
Multiplying (\ref{mamprok3d.continit}) by $-|\boldsymbol{u}_{i}|^2/2$, (\ref{mamprok3d.mominit}) (in the sense of the inner product) by $\boldsymbol{u}_{i}$ and summing, we deduce
$$\frac{\partial}{\partial t}\frac{\rho_i|\boldsymbol{u}_{i}|^2}{2}+{\rm div\,}\left(\frac{\rho_i|\boldsymbol{u}_{i}|^2}{2}\boldsymbol{v}\right)
+{\rm div\,}(\alpha_i p\boldsymbol{u}_{i})-p{\rm div\,}(\alpha_i \boldsymbol{u}_{i})=$$$$={\rm div\,}({\mathbb S}_{i}\boldsymbol{u}_{i})-{\mathbb S}_{i}:{\mathbb D}(\boldsymbol{u}_{i})+
\rho_i\boldsymbol{u}_{i}\cdot\boldsymbol{f}_{i}.$$
Making the summation of these equalities and the equation (\ref{mamprok3d.transpres}) with\linebreak $\varphi(\gamma)=1/(\gamma-1)$, we obtain the desired balance
\begin{equation}\label{mamprok3d.energfin}\begin{array}{c}\displaystyle
\frac{\partial}{\partial t}\left(\sum\limits_{i=1}^N\frac{\rho_i|\boldsymbol{u}_{i}|^2}{2}+\frac{p}{\gamma-1}\right)+
{\rm div\,}\left(\left[\sum\limits_{i=1}^N\frac{\rho_i|\boldsymbol{u}_{i}|^2}{2}+\frac{\gamma p}{\gamma-1}\right]\boldsymbol{v}\right)=\\ \\
\displaystyle ={\rm div\,}\left(\sum\limits_{i=1}^N{\mathbb S}_{i}\boldsymbol{u}_{i}\right)-\sum\limits_{i=1}^N{\mathbb S}_{i}:{\mathbb D}(\boldsymbol{u}_{i})+
\sum\limits_{i=1}^N\rho_i\boldsymbol{u}_{i}\cdot\boldsymbol{f}_{i}.
\end{array}
\end{equation}

\section{Examples of non-monotone pressure}
During further development of the theory, the relation \eqref{mamprok3d.energfin} is insufficient, and the monotonicity of the pressure with respect to the densities is necessary. Let us give examples which show that the equations of the type \eqref{mamprok3d.pressure1} are inconvenient in this sense. To do so, we will stay in the framework of the assumptions that $\gamma>1$ and $K_1=1$ are constants.

\subsection{Non-monotonicity with respect to $\widetilde{\rho}$}
Let $N=2$, $M_1>M_2>0$. Set
\begin{equation}\label{mamprok3d.counterexam1}
\widetilde{\rho}^1=1,\quad \widetilde{\rho}^2=\left(\frac{M_1}{M_2}\right)^{\textstyle\frac{\gamma-1}{2\gamma}}>1,\quad \rho^1=M_1-\varepsilon,\quad
\rho^2=(M_2+\varepsilon)\widetilde{\rho}^2,
\end{equation}
where
$$\varepsilon=\frac{\sqrt{M_1 M_2}}{3}\cdot\frac{M_1-M_2}{M_1+M_2}\in (0,M_1-M_2)\quad\Longrightarrow$$ $$\rho^1\in (M_2,M_1)\quad
\rho^2\in(M_2 \widetilde{\rho}^2,M_1 \widetilde{\rho}^2).$$
Then all components of the corresponding ``reconstructed'' densities
$$\boldsymbol{\rho}=\frac{1}{M_1-M_2}\begin{bmatrix}M_1(\rho-M_2 \widetilde{\rho})\\ M_2(M_1 \widetilde{\rho}-\rho)\end{bmatrix}$$
are positive. However, the corresponding pressures
$$p^1=(\rho^1)^{\gamma-1}\widetilde{\rho}^1=(M_1-\varepsilon)^{\gamma-1},\qquad p^2=(\rho^2)^{\gamma-1}\widetilde{\rho}^2=(M_2+\varepsilon)^{\gamma-1}(\widetilde{\rho}^2)^\gamma,$$
as one can easily see, possess the property $p^2<p^1$, and hence
\begin{equation}\label{mamprok3d.counterexam2}
(p^2-p^1)(\widetilde{\rho}^2-\widetilde{\rho}^1)<0.
\end{equation}

\subsection{Non-monotonicity with respect to $\rho$}
Under the assumptions of the preceding example, we set (\ref{mamprok3d.counterexam1}) again, but let us choose $\varepsilon\in (0,M_1-M_2)$ so that the inequalities
$$\left(\frac{M_1}{M_2}\right)^{\textstyle\frac{\gamma-1}{2\gamma}}<\frac{M_1-\varepsilon}{M_2+\varepsilon}<\left(\frac{M_1}{M_2}\right)^{\textstyle\frac{1}{2}}$$
hold. This gives positive components of the densities again, but now $p^2>p^1$, $\rho^2<\rho^1$, and hence
\begin{equation}\label{mamprok3d.counterexam4}
(p^2-p^1)(\rho^2-\rho^1)<0.
\end{equation}

\subsection{Integral non-monotonicity}
The ``pointwise'' examples above let us construct the corresponding examples of functional density distributions in any domain which satisfy the conservation of the mass, i.~e.
$\int\boldsymbol{\rho}^1=\int\boldsymbol{\rho}^2$, but with the loss of the monotonicity, i.~e. with the integral analogs of the inequalities (\ref{mamprok3d.counterexam2}) and~(\ref{mamprok3d.counterexam4}). To do so, it suffices to split the domain into two parts of equal measure, and then to take the values from the preceding examples on one part of the domain, and to swap densities on the other part. Similarly one can make sure that inequalities of the type
$$\int (p(\boldsymbol{\rho}^1)-p(\boldsymbol{\rho}^2))(\rho^1_j-\rho^2_j)\geqslant 0$$
are impossible, as well as their analogs with any linear combinations of the densities instead of~$\rho_j$.

\subsection{Conclusions concerning admissible constitutive\\ equations}
Hence, if we are based on the constitutive equations of the type (\ref{mamprok3d.pressure1}), i.~e. if the unique (common) pressure depends on all densities but not on their special combination, then we have no hope to obtain inequalities of the form $\overline{\rho p}\geqslant \rho \overline{p}$ (bar stands for the weak limit) which play the crucial role in the EVF technique. In addition, such a situation (one pressure depends on several combinations of the densities) is unlikely for the thermodynamics.

On the other hand, the dependence of the pressures on separate densities would possibly be found reasonable if the momentum equations \eqref{mamprok3d.mominit} contain partial pressures instead of the common one. However, in this case, the present EVF technique works only with diagonal or triangular viscosity matrices, which is of little interest for both physics and mathematics (see details in \cite{mamprok.semr17}).

\section{Further development of the theory}
One can formulate the initial boundary value problem for the system \eqref{mamprok3d.continit},~\eqref{mamprok3d.mominit} in a space-time cylinder, or the steady analog of this problem, and after that, one can construct approximate solution via an appropriate method. Basing on \eqref{mamprok3d.energfin}, after additional efforts, the estimates can be obtained which allow to select weakly convergent sequence of approximate solutions, and the focus is shifted to the justification of the equality $\overline{p}\sim p$ (here $\sim$ means equal action on the test functions) and/or strong convergence of the densities, for which the standard (in the viscous compressible flow theory) EVF technique is applied.

It is important to note that there exist two versions of the final step of the EVF argument. The first one is based on \cite{mamprok.feir09}, P. 339 (it is applied in \cite{mamprok.semi2} and \cite{mamprok.smz2}), and the second one is based on~\cite{mamprok.novstrs04}, P. 188, or on more general arguments as on P. 187 (it is applied in \cite{mamprok.semi1}). In \cite{mamprok.semr17}, a general scheme (in the first version) of applying this technique is written for the systems like \eqref{mamprok3d.continit}, \eqref{mamprok3d.mominit}. This scheme is completely applicable to the system \eqref{mamprok3d.continit}, \eqref{mamprok3d.mominit}, \eqref{mamprok3d.pressure2} provided that all quantities $K$, $\gamma$, $\{\alpha_i\}_{i=1}^N$, $\{\lambda_{ij}\}_{i,j=1}^N$ and $\{\mu_{ij}\}_{i,j=1}^N$ are constant. In the general case, we can state only the following properties of these quantities:
\begin{itemize}
  \item for $K$: \eqref{mamprok3d.transportk},
  \item for $\gamma$: \eqref{mamprok3d.gamma} (possibly, it is sufficient to claim \eqref{mamprok3d.transportgamma}),
  \item for $\{\alpha_i\}_{i=1}^N$: some explicit form of $\alpha_i(\boldsymbol{\xi})$ (possibly, it is sufficient to claim \eqref{mamprok3d.transportalpha}),
  \item for $\{\lambda_{ij}\}_{i,j=1}^N$ and $\{\mu_{ij}\}_{i,j=1}^N$: some explicit form like \eqref{mamprok3d.viscfin} or any other described in Section~3
  (possibly, it is sufficient to claim \eqref{mamprok3d.transportvisc}).
\end{itemize}
Functioning of the EVF technique in the situation when at least one of these quantities is not constant, remains an open problem.

\section{Bibliographic remarks}
As a starting point for the multi-fluid model under consideration, we can cite \cite{mamprok.Nigm} (Chapter 1) and \cite{mamprok.Raj} (Chapter 7). During the development of the existence theory, it became clear that the model admits and requires definite corrections, and we came, as a result, to several variants of the model, one of which is the system \eqref{mamprok3d.continit}, \eqref{mamprok3d.mominit}. The main solvability results for this variant of the model, and for other ones, are obtained in \cite{mamprok.semi2}, \cite{mamprok.smz2}, \cite{mamprok.semi1}, \cite{mamprok.smj12}, \cite{mamprok.izvran}, \cite{mamprok.smz1} and \cite{mamprok.prokkraj}. A detailed review of the results and the discussion of the model can be found in \cite{mamprok.france} and \cite{mamprok.semr17}.

\end{document}